
\documentstyle{amsppt}
\def\L{{\Cal L}}
\def\O{{\Cal O}}

\def\p{{\Bbb P}}

\def\hbar{{\overline h}}

\def\Q{{\Cal Q}}
\def\P{{\Cal P}}
\def\K{{\Cal K}}
\def\M{{\Cal M}}
\def\N{{\Cal N}}
\def\I{{\Cal I}}
\def\E{{\Cal E}}
\def\F{{\Cal F}}
\def\G{{\Cal G}}

\def\Z{{\Bbb Z}}
\def\Pn{{\p}^n}
\def\Pthree{{\p}^3}

\def\ext{\mathop{\Cal Ext}}
\def\hom{\mathop{\Cal Hom}}
\def\ra{\rightarrow}

\topmatter
\title Even Linkage Classes \endtitle
\author Scott R. Nollet \endauthor
\abstract  In this paper the author generalizes the $\E$ and $\N$-type
resolutions used by Martin-Deschamps and Perrin \cite{10} to subschemes of
pure codimension in projective space, and shows that these resolutions
are interchanged by the mapping cone procedure under a simple linkage.
Via these resolutions, Rao's correspondence is extended to give a bijection
between even linkage classes of subschemes of pure codimension two and
stable equivalence classes of reflexive sheaves $\E$ satisfying $H^1_*(\E)=0$
and $\ext^1(\E^\vee,\O)=0$. Further, these resolutions are used to extend
the work of Martin-Deschamps and Perrin for Cohen-Macaulay curves
in $\Pthree$ to subschemes of pure codimension two in $\Pn$. In particular,
even linkage classes of such subschemes satisfy the Lazarsfeld-Rao property
and any minimal subscheme for an even linkage class links directly to a
minimal subscheme for the dual class. \endabstract
\date December 13, 1994 \enddate
\subjclass 13, 14 \endsubjclass
\endtopmatter
\document
\heading Introduction \endheading
   Much progress has been made in the area of linkage theory for subschemes
of $\Pn$ over the past 15 years, but most of the work has dealt
exclusively with locally Cohen-Macaulay (CM) subschemes. This is
somewhat restrictive, since linkage applies more generally to
equidimensional subschemes without embedded components. In particular,
linkage techniques apply to integral subschemes, which need not be
CM. The purpose of this paper is to extend several standard theorems
for CM subschemes to the more general setting. \par
In the first section, we generalize the $\E$ and $\N$-type resolutions
used by Martin-Deschamps and Perrin \cite{10} to subschemes in
$\Pn$ of codimension $r>1$. We show that the equidimensional subschemes
of $\Pn$ having no embedded components can be characterized by an
algebraic condition on their $\E$-type resolutions (see
corollary 1.10). For subschemes which are linked,
we show that the mapping cone procedure \cite{13, proposition 2.5}
interchanges $\E$ and $\N$-type resolutions for the two subschemes. \par
It was a great discovery of Rao that the even linkage classes of locally
Cohen-Macaulay codimension two subschemes of $\Pn$ are in bijective
correspondence with stable equivalence classes of vector bundles $\F$
on $\Pn$ which satisfy $H^1_*(\F)=0$ \cite{15}. In the second section of
this paper, we extend Rao's correspondence to show that the even linkage
classes of pure codimension two subschemes of $\Pn$ are in bijective
correspondence with stable equivalence classes of reflexive sheaves
$\F$ on $\Pn$ which satisfy $H^1_*(\F)=0$ and $\ext^1(\F^\vee,\O)=0$
(see theorem 2.10).
Since $\ext^1(\F^\vee,\O)=0$ for any vector bundle $\F$, this gives
a generalization of vector bundle. \par
Perhaps the most important theorem about even linkage classes of locally
Cohen-Macaulay codimension two subschemes of $\Pn$ is the fact that
they satisfy the Lazarsfeld-Rao property (see \cite{1}; the case of
curves in $\Pthree$ was done in \cite{10}). By replacing vector bundles
with reflexive sheaves, the proof in \cite{10} goes through with
minor modifications (see theorem 3.26). The concept of domination
is used to simplify the statement of the main theorem.
\par
The last two sections of this paper deal only with subschemes of codimension
two.  It would be interesting to know if
Rao's correspondence (section two) can be generalized in some way via
the resolutions used in section one. If this situation were better
understood, then perhaps the methods of section four might be used
to investigate whether the Lazarsfeld-Rao property holds for even
linkage classes of subschemes of codimension $>2$. \par
I'd like to thank Robin Hartshorne and Juan Migliore for editorial
suggestions, especially regarding the first section. I appreciate
personal communications with A. P. Rao about the details of
his work. Finally, I'd like to thank Martin-Deschamps and Perrin for
their book \cite{10}. The generality with which they handle the
fourth chapter saved this author a great deal of work. \par
%
%
\heading Two Cone Constructions \endheading
In this section we generalize the $\E$-type and $\N$-type resolutions
for curves found in \cite{10} to subschemes of pure codimension $r$ in
$\Pn$ (in other words, equidimensional subschemes with no embedded
components). We show in this broader setting
that the mapping cone procedure \cite{13, proposition 2.5} still works and
interchanges these two types of resolutions for subschemes which are linked.
Before getting on to these results we briefly recall the definition of
linkage and some of the first properties. The main reference for this is
\cite{7}, where the foundations of linkage theory have been rewritten in
the more conceptual context of generalized divisors. \par
\definition{Definition 1.1} Let $V_1,V_2$ be subschemes of codimension
$r$ in $\Pn$.
We say that $V_1$ {\it is linked to} $V_2$ {\it by the complete
intersection} $X$ if $X$ is a global complete intersection of codimension
$r$ containing $V_1$ and $V_2$ in such a way that
$$\I_{V_2,X} \cong \hom(\O_{V_1},\O_X)$$
$$\I_{V_1,X} \cong \hom(\O_{V_2},\O_X)$$
\enddefinition
\proclaim{Proposition 1.2} Linkage enjoys the following properties:
\roster
\item If $V_1$ is linked to $V_2$ by the complete intersection $X$,
then $V_2$ is linked to $V_1$ by $X$.
\item If $V_1 \subset \Pn$ is of pure codimension $r$ and is contained in a
complete intersection $X$ of codimension $r$, then there exists a unique
subscheme $V_2 \subset \Pn$ of pure codimension $r$ such that $V_2$ is
linked to $V_1$ by $X$.
\item If $V_1$ is linked to $V_2$ by $X$, then either both $V_1$ and $V_2$ are
of
pure codimension $r$ in $\Pn$, or one of them is empty.
\item If $V_1$ is linked to $V_2$, then $V_1$ is locally Cohen-Macaulay if
and only if $V_2$ is locally Cohen-Macaulay.
\endroster
\endproclaim
\demo{Proof}
Part (3) is proved in \cite{12, proposition 4.2 (2)}, and
credited to Schwartau's thesis, although it was surely known before then.
Parts (1), (2) and (4) can be found in propositions 4.1 and 4.2 of \cite{7}.
In the latter paper, the definition of linkage is different, but is proved
to be equivalent to the usual definition.
\enddemo
\remark{Remark 1.3}
Part (4) shows that linkage preserves the property of
being locally Cohen-Macaulay. In the literature, many papers only
apply linkage theory to locally Cohen-Macaulay schemes. The purpose of
this paper is to generalize many of these results to subschemes which are
not locally Cohen-Macaulay.
\endremark
\definition{Definition 1.4}
A sheaf $\F$ on $\Pn$ is called {\it dissoci\'e} if
it is a direct sum of line bundles.
\enddefinition
\definition{Definition 1.5}
If $V \subset \Pn$ is a subscheme of
codimension $r$, then an {\it $\E$-type resolution for} $X$ is an exact
sequence
$$0 \ra \E \ra \Q_{r-1} \ra \Q_{r-2} \ra \dots \ra \Q_1 \ra \O
(\ra \O_V \ra 0)$$
such that $\Q_i$ is dissoci\'e for $0 < i < r$ and $H^i_*(\E)=0$ for
$0 < i < r$.
\enddefinition
\remark{Remark 1.6}
Note that any closed subscheme $V \subset \Pn$ has an $\E$-type
resolution. This is most easily constructed by sheafifying a graded free
resolution of $I_V$, and letting $\E$ be the kernel at the $r$th map.
If $F_i$ is the kernel of the map from $\Q_i$, then breaking the
resolution into short exact sequence shows that for each $0 \leq k < r$,
the graded modules $H^i_*(F_{i+k})$ are all isomorphic for $1 \leq i < r-k$.
Thus the condition $H^i_*(\E)=0$ for $0 < i < r$ is equivalent to the
condition that $H^1_*(F_i)=0$ for $1 \leq i < r$, which is equivalent to
saying that applying $H^0_*$ to the resolution gives a resolution of
graded $S$-modules for $I_V$.
\endremark
\remark{Remark 1.7}
Since $V$ is a scheme of dimension $n-r$, the local depth of $\O_V$ is
$\leq n-r$, hence the local cohomological dimension of $\O_V$ is always
$\geq r$ by the Auslander-Buchsbaum theorem. This shows that $\E \neq 0$.
Also, $\E$ is locally free if and only if $V$ is locally
Cohen-Macaulay, because both these conditions are equivalent to the
condition that the local cohomological dimension of $V$ is equal to $r$.
\endremark
\proclaim{Proposition 1.8}
Let $V \subset \Pn$ be a subscheme which is contained
in a complete intersection $X$ of $r$ hypersurfaces
of degrees $\{d_i\}$. Let $G$ be an $\E$-type resolution for $V$ and let
$F$ denote the Koszul resolution for $X$. Then there is a morphism of
complexes $\alpha:F \ra G$ induced by the inclusion $\I_X \ra \I_V$, and
the mapping cone of $\alpha^\vee(-\sum{d_i})$ gives a resolution for $\O_W$,
where $W$ is the subscheme of $X$ defined by $\I_{W,X} \cong \hom(\O_V,\O_X)$.
\endproclaim
\demo{Proof}
In the context of definition 1.5, we identify
$G_0$ with $\O$ and $G_r$ with $\E$. By the observation made in
remark 1.6, the condition $H^i_*(\E)=0$ for $0 < i < r$ is
equivalent to the condition that applying $H^0_*$ to the resolution gives
an exact sequence of graded $S$-modules. In particular, we can find a morphism
from the Koszul resolution of graded $S$-modules to the $\E$-type resolution,
and sheafifying gives the desired morphism $\alpha$.
\par
Now let $d=\sum{d_i}$. Letting $C$ be the mapping cone of $\alpha^\vee(-d)$,
we get a long exact homology sequence for the three complexes involved.
{}From \cite{5, III, proposition 6.5} we know that $\ext$ can be computed from
locally free resolutions. It follows that $H_i(F^\vee(-d))=\ext^i(\O_X,\O)$.
Since $X$ is of codimension $r$, these $\ext$s vanish for $i < r$
\cite{5, III, lemma 7.1}, and hence the long exact homology sequence
gives isomorphisms $H_i(C) \cong H_i(G^\vee(-d)$ for $0 \leq i <r$ and an
exact sequence
$$0 \ra H_r(C) \ra H_r(G^\vee(-d)) \ra H_r(F^\vee(-d)) \ra H_{r+1}(C) \ra 0$$
\par
Now let $L$ be a dissoci\'e resolution for $\O_V$ which agrees with $G$
up to the $(r-1)$st step, so that $L_{\geq r}$ is a resolution for $G_r=\E$,
$L_0=\O$ and $L_i=\Q_i$ for $0 < i < r$. Clearly the morphism $\alpha$ can
be factored through the complex $L$ as
$F @>\gamma>> L @>\beta>> G$. Since $\beta$ is a
morphism of locally free resolutions of $\O_V$ and $\O_X$, the induced homology
map $H_r(\beta^\vee(-d))$ is the canonical map
$\ext^r(\O_V,O)(-d) \ra \ext^r(\O_X,O)(-d)$. Since the functor
$\ext^r(*,\O)(-d)$ on the category of $\O_X$-modules is isomorphic to the
functor $\hom(*,\O_X)$, this map is identified with the inclusion
$\hom(\O_V,\O_X) \hookrightarrow \hom(\O_X,\O_X)$ and it follows that the map
$H_r(L^\vee(-d)) \ra H_r(F^\vee(-d))$ is injective and has $\O_W$ as cokernel.
Thus, to finish the proof of the proposition, it suffices to show that the
homology maps $H_i(G^\vee(-d)) \ra H_i(L^\vee(-d))$ are isomorphisms for
$0 \leq i \leq r$. \par
For $0 \leq i < r-1$ it is clear that $H_i(\beta^\vee(-d))$ is an isomorphism,
because $\beta$ is an isomorphism in these degrees. For $i=r-1$ and $r$, we
have a fragment of the commutative diagram corresponding to $\beta^\vee$
$$ \CD
\Q_{r-2}^\vee @>>> \Q_{r-1}^\vee @>>> \E^\vee @>>> 0 \\
@VVV  @VVV  @VVV  @VVV \\
L_{r-2}^\vee @>>> L_{r-1}^\vee @>>> L_r^\vee @>>> L_{r+1}^\vee
\endCD $$
By our choice of $L$ we have that the two vertical maps on the left are
isomorphisms. Because $L_{\geq r}$ is a resolution for $\E$ and dualizing is
left exact, we see that $\E^\vee$ is precisely the kernel of the map
$L_r^\vee \ra L_{r+1}^\vee$. It follows from commutativity of the diagram
that the homology maps $H(\Q_{r-1}^\vee) \ra H(L_{r-1}^\vee)$ and
$H(\E) \ra H(L^\vee)$ are isomorphisms. \par
\enddemo
\remark{Remark 1.9}
If we skip past the first paragraph, the proof above
is local commutative algebra. The local statement proved is as follows:
Let $(A,m)$ be a regular local ring and let $I \subset A$ be an ideal
which contains a regular sequence $\{a_i\}$ of length $r$. Let $J$
be the ideal generated by the $a_i$. If $C$ is a resolution for $A/I$ of
the form
$$0 \ra E \ra F_{r-1} \ra F_{r-2} \ra \dots \ra F_1 \ra A (\ra A/I \ra 0)$$
with $F_i$ free, $D$ is the Koszul complex for $A/J$, and $\alpha$
is any morphism $D \ra C$ induced by the map $A/J \ra A/I$, then
the mapping cone of $\alpha^\vee$ gives a resolution for $A/K$, where
$K$ is the ideal containing $J$ defined by $K/J \cong \hom(A/I,A/J)$.
\endremark
\proclaim{Corollary 1.10}
If $V \subset \Pn$ is a subscheme of pure
codimension $r$ with $\E$-type resolution as in definition 1.5
and $V$ links to $W$ via a complete intersection $X$ of hypersurfaces with
degrees summing to $d$ and having Koszul resolution $F$, then $W$ has a
resolution
$$0 \ra \P_{r-1} \ra \P_{r-2} \ra \dots \ra \P_1 \ra \N \ra \O \ra \O_W \ra 0$$
where $\P_{r-1} = \Q_1^\vee(-d)$,
$\P_i=\Q_{r-i}^\vee(-d) \oplus F_{r-i-1}^\vee(-d)$ for $1 \leq i < r-1$ and $\N
= \E^\vee(-d) \oplus F_{r-1}^\vee(-d)$.
\endproclaim
\demo{Proof} In this case, the scheme $W$ linked to $V$ satisfies
$\I_{W,X} \cong \hom(\O_V,\O_X)$ by the definition of linkage.
Hence it is the scheme $W$ of the previous proposition.
The mapping cone obtained in the previous proposition has form
almost like the one given above, except that the summands
$\O^\vee(-d)$ at the left end have been split off (which is possible
because the induced map is the identity).
\enddemo
\proclaim{Proposition 1.11}
Let $V \subset \Pn$ be a subscheme which is contained in the
complete intersection $X$ of $r$ hypersurfaces of degrees $\{d_i\}$.
Let $G$ be a resolution for $\O_V$ of the form
$$0 \ra \P_{r-1} \ra \P_{r-2} \ra \dots \ra \P_1 \ra \N \ra \O
(\ra \O_V \ra 0)$$
where $\P_i$ are dissoci\'e. If $F$ denotes the Koszul resolution for $X$,
then there is a morphism $\alpha:F \ra G$ induced by the inclusion
$\I_X \ra \I_V$. If $W$ is the subscheme of $X$ defined by
$\I_{W,X} \cong \hom(\O_V,\O_X)$, then the mapping cone of the
morphism $\alpha^\vee(-\sum{d_i})$ gives a resolution of $\O_W$ if and
only if $\ext^i(\N,\O)=0$ for $1 \leq i < r$.
\endproclaim
\demo{Proof} Let $R_i$ denote the image of $\P_i$ under the maps for the
resolution for $\O_V$. We can break the resolution for $\O_V$ into short
exact sequences (with dissoci\'e sheaves in the middle). The long
exact cohomology sequence shows that for each
$0 < i <r$ we have that $H^1_*(R_i) \cong H^{r-i}_*(R_{r-1})$, but
$R_{r-1} \cong \P_{r-1}$ is dissoci\'e, hence all these graded
$S$-modules are zero. It follows that when we apply $H^0_*$ to $G$
we get an exact sequence of graded $S$-modules. We can find a resolution
for $I_V$ consisting of free graded $S$-modules which maps surjectively
to this exact sequence on each component. Sheafifying gives a dissoci\'e
resolution $L$ for $\O_V$ which maps surjectively to the resolution that
we started with. The inclusion of total ideals $I_X \subset I_V$ induces
a morphism from the Koszul resolution of graded $S$-modules for
$\O_X$ to the graded free resolution $H^0_*(L)$. Thus we obtain a morphism
$\alpha:F \ra G$ of complexes which factors through $L$ as
$$F @>\gamma>> L @>\beta>> G.$$ \par
Let $d=\sum{d_i}$. Letting $C$ be the mapping cone of $\alpha^\vee(-d)$,
we obtain a long exact homology sequence. As in the proposition 1.8
the fact that $V$ is of codimension $r$ shows that we have isomorphisms
$H_i(C) \cong H_i(G^\vee(-d))$ for $0 \leq i < r$ and an exact sequence
$$0 \ra H_r(C) \ra H_r(G^\vee(-d)) \ra H_r(F^\vee(-d)) \ra H_{r+1}(C) \ra 0$$
Computing $\ext$ from the locally free resolution $L$ gives
$$H_r(L^\vee(-d)) = {\ext}^r(\O_V,\O)(-d)=\hom(\O_V,\O_X)$$
and this last sheaf is isomorphic to $\I_{W,X}$ by definition of
$W$. The map from this sheaf to $\O_X$ is the canonical inclusion
because the functor ${\ext}^r(*,\O)(-d)$ is naturally isomorphic to
the functor $\hom(*,\O_X)$ on the category of $\O_X$-modules.
It follows that the mapping cone to $\alpha^\vee(-d)$ is a resolution
for $\O_W$ if and only if the homology map $H^i(\beta^\vee)$ is an
isomorphism for $0 \leq i \leq r$.
\par
Now we compare the homology of $L^\vee$ and $G^\vee$. Let $K$ be the
kernel of the complex map $\beta$. Since we chose $L$ with $L_0=\O$,
we have $K_0=0$.
Now dualize the exact sequence of complexes to get a diagram
$$ \CD
0 @>>> G_0^\vee @>>> G_1^\vee @>>> G_2^\vee @>>> \dots \\
@. @VVV  @VVV  @VVV \\
0 @>>> L_0^\vee @>>> L_1^\vee @>>> L_2^\vee @>>> \dots \\
@. @VVV  @VVV  @VVV \\
@. 0 @>>> R @>>> K_2^\vee @>>> \dots
\endCD $$
where $R$ denotes the image of the morphism $L_1^\vee \ra K_1^\vee$.
For $i>1$, $K_i$ is the kernel of the surjection $L_i \ra G_i$ of bundles,
hence the $K_i$ are locally free and the $i$th column of the
diagram is a short exact sequence. Clearly the first two columns are
exact (the second is exact by the choice of $R$), so in fact all the
columns are exact and we have an associated long exact homology sequence.
Since $G_i=0$ (and hence $H_i(G^\vee)=0$) for $i>r$, this long exact sequence
shows that $H_i(\beta^\vee)$ is an isomorphism for
$0 \leq i \leq r$ if and only if the homology $H_i$ of the bottom
complex is zero for $i \leq r$. \par
Noting that $R$ is a submodule of $K_1^\vee$, $H(R)=0$ in any event.
Since dualizing is left exact, we see that $H(K_2^\vee)=0 \iff R=K_1^\vee$.
By definition of $R$, this happens $\iff L_1^\vee \ra K_1^\vee$ is
surjective, but the cokernel of this map is just $\ext^1(\N,\O)$,
so $H(K_2^\vee)=0 \iff \ext^1(\N,\O)=0$. Finally, since $K_1$ is the
kernel of the map $L_1 \ra \N$ we see that
$$\dots \ra K_3 \ra K_2 \ra L_1$$
is a locally free resolution for $\N$. Since $\ext$ can be computed from
such resolutions, we see that $H_i(K^\vee) \cong \ext^{i-1}(\N)$ for $i>1$.
It follows that $H_i(K^\vee)=0$ for $2 < i \leq r$ if and only if
$\ext^i(\N)=0$ for $1 < i \leq r-1$. Combining the conclusions of this
paragraph shows that the homology $H_i$ of the bottom complex is exact
for $i \leq r \iff \ext^i(\N)=0$ for $1 \leq i < r$. This proves the
proposition.
\enddemo
\remark{Remark 1.12}
Again we have proved a local algebra result: Let $(A,m)$ be
a regular local ring, and let $I$ be an ideal which contains a regular
sequence $\{a_i\}$ of length $r$ which generates the ideal $J$. Let $K$
be the ideal containing $J$ defined by the condition $K/J \cong \hom(A/I,A/J)$,
and let $C$ be a resolution
$$0 \ra F_{r-1} \ra F_{r-2} \ra \dots \ra F_1 \ra N \ra A (\ra A/I \ra 0)$$
for $A/I$ with $F_i$ free. If $\alpha$ is a morphism from the Koszul
resolution of $A/J$ to $C$ induced by the surjection $A/J \ra A/I$, then
the mapping cone of $\alpha^\vee$ gives a resolution for $A/K$ if and only
if $\ext^i(N,A)=0$ for $0 < i < r$.
\endremark
\proclaim{Proposition 1.13}
Let $V \subset \Pn$ be of pure codimension $r$ with
$\E$-type resolution
$$0 @>>> \E @>>> \Q_{r-1} @>>> \dots @>>> \Q_1 @>\pi>> \O
(@>>> \O_V @>>> 0). $$
Then $\ext^i(\E^\vee,\O)=0$ for $0 < i < r$.
\endproclaim
\demo{Proof} The idea here is simple. We first use a complete intersection
$X$ to link $V$ to another scheme $W$, and use proposition 1.8
to get a resolution for $W$. Then we use a carefully chosen morphism of
complexes
from the Koszul complex for $\O_X$ to this resolution, and compute
the cone of the dual of this morphism. By the choice of morphism, we can
cancel some redundant factors and are left with the $\E$-type resolution
that we started with.
Applying proposition 1.11 to the resolution for $\O_W$ shows that
$\ext^i(\E^\vee,\O)=0$ for $0 < i < r$. \par
Proceeding as suggested above, let $X$ be a complete intersection of $r$
hypersurfaces of degrees $d_i$ which links $V$ to another subscheme $W$.
Let $F$ denote the Koszul complex for $\O_X$, so that
$F_0=\O, F_1=\oplus\O(-d_i)$ and $F_r=\O(-\sum{d_i})$.
By looking at the maps involved, it's easy to see that the dual
complex $F^\vee$ is a Koszul resolution for $\O_X \otimes F_r^\vee$.
Let $\alpha$ denote the map of complexes from the Koszul resolution $F$
of $\O_X$ to the $\E$-type resolution of $\O_V$.
Dualizing $\alpha$ and applying the mapping cone procedure (proposition
1.8) gives a resolution for $\O_W \otimes F_r^\vee$.
This complex, along with a morphism of the (dual) Koszul complex for
$\O_X \otimes F_r^\vee$ to this complex, is given in the
diagram below.
$$ \CD
\O^\vee @>>> \Q_1^\vee \oplus F_0^\vee @>>> \Q_2^\vee \oplus F_1^\vee @>>>
\dots @>>> F_r^\vee \\
@. @AAA @AAA @AAA @AAA \\
@. F_0^\vee @>>> F_1^\vee @>>> \dots @>>> F_r^\vee
\endCD $$
The differential on the mapping cone complex is given by the usual formula,
namely, $\delta(x,y)=(\delta(x),\alpha(x)-\delta(y))$. Given that this is
the differential, it is easy to see that the diagram commutes if we take the
morphism of complexes to be the second projection. Further, we can cancel
off the redundant copy of $\O^\vee$ and $F_0^\vee$ on the left if we replace
the
vertical map from $F_0^\vee$ to $\Q_1^\vee \oplus F_0^\vee$ with the map
$F_0^\vee=\O^\vee @>\pi^\vee>> \Q_1^\vee$. \par
Using this morphism of complexes, we apply proposition 1.11.
We dualize the modified morphism of complexes. After removing
all the double duals, we get a diagram
$$ \CD
F_r @>>> \E \oplus F_{r-1} @>>> \dots @>>> \Q_1 \\
@VVV         @VVV               @VVV       @VVV \\
F_r @>>> F_{r-1}           @>>> \dots @>>> F_0
\endCD $$
Here the leftmost vertical map is the identity, the middle vertical maps
are projections from the second factor, and the rightmost vertical map
is $\pi$. After cancelling off the extra summand of $F_r$ that appears,
the mapping cone takes the form
$$0 \ra \E \oplus F_{r-1} \ra \Q_{r-1} \oplus F_{r-2} \oplus F_{r-1} \ra
\dots \Q_1 \oplus F_1 \ra F_0$$
The induced map on the summands $F_i$ is the identity, so these can also
be removed. Looking at the induced maps, we see that this is the complex
for $\O_X$ that we started with. Since $V$ is linked to $W$, we can apply
proposition 1.11 to see that
$\ext^i((\E^\vee \oplus F_{r-1}^\vee)\otimes F_r,\O)=0$ for $0 < i < r$.
Since $F_{r-1}^\vee$ and $F_r$ are dissoci\'e, we conclude that
$\ext^i(\E^\vee,\O)=0$ for $0< i < r$. This finishes the proof.
\enddemo
\remark{Remark 1.14}
Put into the context of commutative algebra, the proof above
yields the following result: Let $A$ be a regular local ring, $I \subset A$ an
ideal such that all the associated primes to $A/I$ have height $r$.
If $$0 \ra E \ra F_{r-1} \ra F_{r-2} \ra \dots \ra F_0 (\ra A/I \ra 0)$$
is a resolution of $A/I$ with $F_i$ free, then $Ext^i(Hom(E,A),A)=0$ for each
$0 < i < r$.
\endremark
\definition{Definition 1.15} If $V \subset \Pn$ is a subscheme of
codimension two, then an {\it $\N$-type resolution for} $\O_V$ is an
exact sequence
$$0 \ra \P_{r-1} \ra \P_{r-2} \ra \dots \P_1 \ra \N \ra \O (\ra \O_V \ra 0)$$
where $\P_i$ are dissoci\'e, $\N$ is reflexive, and $\ext^i(\N,\O)=0$ for
$0<i<r$.
\enddefinition
\proclaim{Corollary 1.16}
Let $V \subset \Pn$ be a nonempty subscheme of
codimension $\geq r$. Then $V$ is of pure codimension $r$ if and only if
$V$ has an $\E$-type resolution
$$0 \ra \E \ra \Q_{r-1} \ra \Q_{r-2} \ra \dots \Q_1 \ra \O (\ra \O_V \ra 0)$$
such that $\ext^i(\E^\vee,\O)=0$ for $0<i<r$.
\endproclaim
\demo{Proof} The previous proposition gives us the forward direction.
For the converse, suppose that $V$ has such a resolution. Since $V$ has
codimension $\geq r$, we can find a complete intersection scheme $X$
which contains $V$. Letting $\alpha$ be
a morphism from the Koszul resolution of $\O_X$ to the
$\E$-type resolution for $X$, we apply proposition 1.8 to
see that the mapping cone of $\alpha^\vee(-d)$ gives a resolution for
$W \subset X$ defined by
$\I_{W,X} \cong \hom(\O_V,\O_X)$. By our condition on
$\E$, we see that this is an $\N$-type resolution for $W$. \par
Now we can apply proposition 1.11 to the resolution for $W$.
If we use the morphism of complexes chosen in the proof of proposition
 1.13, we get an $\E$-type resolution for a scheme $V^\prime$
defined by $\I_{V^\prime,X} \cong \hom(\O_W,\O_X)$. By choice of
the morphism, we see that the resolution can be simplified
to agree with our starting resolution for $V$. It follows that $V^\prime = V$
and that $V$ and $W$ are linked by the complete intersection $X$.
In particular, $V$ is pure codimension $r$.
\enddemo
\remark{Remark 1.17}
Again we give the corresponding commutative algebra
statement: Let $(A,m)$ be a regular local ring, and let $I \subset A$ be
an ideal such that all the associated primes of $A/I$ have height $\geq r$.
If $$0 \ra E \ra F_{r-1} \ra F_{r-2} \ra \dots \ra F_0 (\ra A/I \ra 0)$$
is a resolution of $A/I$ with $F_i$ free, then all associated primes
of $A/I$ have height exactly $r \iff Ext^i(Hom(E,A),A)=0$ for $0 < i < r$.
This local version shows that the global version can be strengthened in
the following way: if we only assume that the $\Q_i$ are locally free,
then $V$ has pure codimension $r \iff \ext^i(\E^\vee,\O)=0$ for $0 < i < r$.
\endremark
\proclaim{Corollary 1.18}
Let $V \subset \Pn$ be a nonempty subscheme of
codimension $\geq r$. Then $V$ is of pure codimension $r$ if and only if
$V$ has an $\N$-type resolution
$$0 \ra \P_{r-1} \ra \P_{r-2} \ra \dots \P_1 \ra \N \ra \O (\ra \O_V \ra 0)$$
such that $H^i_*(\N^\vee)=0$ for $0<i<r$.
\endproclaim
\demo{Proof} If $V$ is of pure codimension $r$, we can link $V$
to a subscheme $W$ of pure codimension $r$. If we use an $\E$-type
resolution for $W$ and apply propositions 1.13 and 1.11, we
get an $\N$-type resolution for $V$ as described in the corollary.
\par
Conversely, suppose that $V \subset \Pn$ has an $\N$-type resolution
as described. Since $V$ has codimension $\geq r$, we can find a complete
intersection $X$ which contains $V$. Let $\alpha$ be a morphism
from the Koszul complex of $\O_X$ to the $\N$-type resolution
of $V$. Applying proposition 1.11, we see that the mapping cone
of $\alpha^\vee(-d)$ gives an $\E$-type resolution for a subscheme
$W \subset X$ defined by
$\I_{W,X} \cong \hom(\O_V,\O_X)$. As in the proof of
proposition 1.13, we can find a morphism from the Koszul
complex of $\O_X$ to this $\E$-type resolution in such a
way that the cone of the dual morphism twisted by $-d$ simplifies
to the original $\N$ type resolution for $V$. It follows (as in
the previous corollary) that $V$ is linked to $W$, hence $V$ is
of pure codimension $r$.
\enddemo
\remark{Remark 1.19}
In comparing the last two corollaries with the situation
when $V$ is locally Cohen-Macaulay, and hence has locally free
resolutions of the two forms suggested, we see that reflexive sheaves
$\E$ such that $\ext^i(\E^\vee,\O)=0$ for $0 < i < r$ give a nice
generalization of vector bundles, as do their duals, the reflexive sheaves
$\N$ such that $\ext^i(\N,\O)=0$ for $0 < i < r$.
\endremark
%
%
\heading Rao's Correspondence \endheading
In this section, we extend Rao's correspondence between stable equivalence
classes of vector bundles and even linkage classes of Cohen-Macaulay
subschemes of codimension two. By specializing the results of the previous
section to the case $r=2$, we find that we can remove the Cohen-Macaulay
hypothesis on the even linkage classes if we use stable equivalence classes
of reflexive sheaves $\E$ such that $H^1_*(\E)=0$ and $\ext^1(\E^\vee,\O)=0$.
Since only minor modifications are necessary,
we keep the outline of Rao's proof \cite{15} with few changes. \par
\definition{Definition 2.1}
The equivalence relation generated by linking is
called {\it liaison} or {\it linkage}. If two schemes can be connected
by an even number of linkages, we say that they are {\it evenly linked}
and the corresponding equivalence relation is called {\it biliaison} or
{\it even linkage}.
\enddefinition
\definition{Definition 2.2} Two reflexive sheaves $\E_1$ and $\E_2$ on $\Pn$
are {\it stably equivalent} if there exist dissoci\'e sheaves $\Q_1,\Q_2$
and $h \in \Z$ such that $\E_1 \oplus \Q_1 \cong \E_2(h) \oplus \Q_2$.
This is an equivalence relation among reflexive sheaves on $\Pn$.
\enddefinition
\proclaim{Proposition 2.3}
Let $\F$ be a reflexive sheaf on $\Pn$, and let $\Omega$ denote the
stable equivalence class determined by $\F$. Let $\F_0$ be a sheaf in
$\Omega$ of minimal rank. Then for each $\G \in \Omega$, there exists
$h \in \Z$ and a dissoci\'e sheaf $\Q$ such that $\G \cong \F_0(h) \oplus \Q$.
In particular, the elements of minimal rank are unique up to twist.
\endproclaim
\demo{Proof} Let $\F_0$ be as in the statement, and let $\G \in \Omega$
be arbitrary. Since $\F_0$ and $\G$ are in the same stable equivalence
class, there exists $h \in \Z$ and dissoci\'e sheaves $\Q_1,\Q_2$ such
that $\F_0(h) \oplus \Q_1 \cong \G \oplus \Q_2$. Let $\O(-a)$ be a
line bundle summand of $\Q_2$ and write $\Q_2 = \O(-a) \oplus \Q_2^\prime$
with $\Q_2^\prime$ dissoci\'e. If we twist the isomorphism by $a$ and
consider the map on global sections, we obtain a surjection
$$H^0(\F_0(h+a)) \oplus H^0(\Q_1(a)) \ra H^0(\O)=k$$
If the induced map $H^0(\F_0(h+a)) \ra k$ is surjective, then the corresponding
map $\F_0(h+a) \ra \O$ is a split surjection. If $\E$ is the kernel, then
$\E \oplus \O \cong \F_0(h+a)$ shows that $\F_0$ was not of minimal rank in
$\Omega$, a contradiction. It follows that the induced map zero, and
hence the map $H^0(\Q_1(a)) \ra k$ is surjective. This shows that the
induced map $\Q_1 \ra \O(-a)$ is split surjective, and if $\Q_1^\prime$ is
the kernel, then we can split off $\O(-a)$ from the original isomorphism to get
$$\F_0(h) \oplus \Q_1^\prime \cong \G \oplus \Q_2^\prime$$
Continuing in this fashion, we can split off all the summands of $\Q_2$ to
obtain an isomorphism $\G \cong \F_0(h) \oplus \Q$ as in the proposition.
\enddemo
\proclaim{Proposition 2.4} There is a well-defined map $\Phi$ from
even linkage classes of pure codimension two subschemes of $\Pn$ to
stable equivalence classes of reflexive sheaves $\E$ such that $H^1_*(\E)=0$
and $\ext^1(\E^\vee,\O)=0$. This map is given by taking $X$ to the last sheaf
in an $\E$-type resolution for $X$.
\endproclaim
\demo{Proof} The map has been given above. It is clear from the definition
of $\E$-type resolution and proposition 1.13 that the image of
$\Phi$ consists of reflexive sheaves $\E$ satisfying $H^1_*(\E)=0$ and
$\ext^1(\E^\vee,\O)=0$. We need only show that it is
well-defined. For this, assume that $X_1,X_2$ are in the same even linkage
class. We can get from one to the other by an even number of links.
If $X_1$ has an $\E$-type resolution
$$0 \ra \E_1 \ra \Q_1 \ra \O (\ra \O_{X_1} \ra 0)$$
then by applying the cone construction an even number of times, we get
a resolution for $X_2$ of the form
$$0 \ra \E_1(h) \oplus \P_2 \ra \Q_2 \ra \O (\ra \O_{X_2} \ra 0)$$
where $\P_2$ and $\Q_2$ are dissoci\'e. Clearly $\E_1$ and
$\E_1(h) \oplus \P_2$ are in the same stable equivalence class, so it
suffices to show that if $X_2$ has a second $\E$-type resolution
$$0 \ra \F \ra \Q \ra \I_{X_2} \ra 0$$
then $\F$ is stably equivalent to $\E_1(h) \oplus \P_2$. \par
For this, we add the map $\Q_2 \ra \I_{X_2}$ to the resolution
involving $\F$, we obtain a commutative diagram
$$\CD
0 @>>> \F @>>> \Q @>>> \I_{X_2} @>>> 0 \\
@. @VVV @VVV @VVV @. \\
0 @>>> \K @>>> \Q \oplus \Q_2 @>>> \I_{X_2} @>>> 0 \\
@. @VVV @VVV @. @. \\
@. \Q_2 @>>> \Q_2 @. @.
\endCD $$
where the columns and rows are exact. Since $H^1_*(\F)=0$ by assumption,
the map $H^0_*(\K) \ra H^0_*(\Q_2)$ is surjective. Since $\Q_2$ is
dissoci\'e, it follows that the map $\K \ra \Q_2$ is split surjective,
whence $\K \cong \F \oplus \Q_2$. Applying the same argument to the
other resolution of $\I_{X_2}$ shows that $\K \cong \E_1(h) \oplus \P_2 \oplus
\Q$, so that
$$\F \oplus \Q_2 \cong \E_1(h) \oplus \P_2 \oplus \Q$$
which finishes the proof.
\enddemo
\proclaim{Lemma 2.5}
Let $X \subset \Pn$ be a subscheme of pure codimension
two having an $\N$-type resolution
$$0 \ra \P \ra \N \ra \I_X \ra 0$$
with $H^1_*(\N^\vee)=0$. If $\N \cong \N^\prime \oplus \O(-p)$, then there
exists a subscheme $X^\prime$ in the even linkage class of $X$ which has a an
$\N$-type resolution of the form
$$0 \ra \P^\prime \ra \N^\prime(d) \ra \I_{X^\prime} \ra 0$$
for some $d \in \Z$.
\endproclaim
\demo{Proof} First note that $\ext^1(\N^\prime,\O)=0$, because $\ext$ respects
direct sums. Moreover, the fact that dualizing respects direct sums shows
that $\N^\prime$ is reflexive. These two facts show that it is at least
possible
to have an $\N$-type resolution involving $\N^\prime$. \par
Let $s \in H^0(\N(p))$ be a nowhere vanishing section
corresponding to the summand $\O(-p)$ of $\N$ and let $\sigma:\N(p) \ra \O$ be
a fixed splitting of $s$ (The quotient of $\N(p)$ by $s$ is $\N^\prime(p)$).
Let $f$ be the image of $s$ in $I_X$. As in \cite{14, proposition 1.11}, we
consider two cases. \par
If $f=0$, then the map $\pi$ factors through $\N^\prime$, let
$\pi^\prime:\N^\prime \ra \I_X$ be the induced surjection. Let $\P^\prime$
be the kernel of $\pi^\prime$. This gives a commutative diagram
$$\CD
0 @>>> \P^\prime @>>> \N^\prime @>>> \I_X @>>> 0 \\
@. @VVV @VVV @VVV @. \\
0 @>>> \P @>>> \N^\prime \oplus \O(-p) @>>> \I_X @>>> 0 \\
@. @VVV @VVV @VVV @. \\
0 @>>> \O(-p) @>>> O(-p) @>>> 0 @.
\endCD $$
Since the map of graded modules $H^0_*(\N^\prime) \oplus H^0_*(\O(-p)) \ra I_X$
is surjective and the induced
map $H^0_*(\O(-p)) \ra I_X$ is zero, the induced map
$H^0_*(\N^\prime) \ra I_X$ is also surjective. Since we also have an
isomorphism
$H^1_*(\N^\prime) \cong H^1_*(\I_X)$ (both are isomorphic to $H^1_*(\N)$), we
see that $H^1_*(\P^\prime)=0$. The snake lemma shows that the vertical
exact sequence on the left is exact. Since $\P$ is dissoci\'e, this sequence
shows that $H^i_*(\P^\prime)=0$ for $2 \leq i < n$. Combining with the
previous observation, we see that all the middle cohomology of $\P^\prime$
vanishes. Applying \cite{4, lemma 6.3}, we conclude that $\P^\prime$
is dissoci\'e. This finishes the case $f=0$. \par
If $f \neq 0$, then it gives the equation of a hypersurface $S$ of degree
$s$ which contains $X$. Let $T$ be another hypersurface
(of degree $t$) which meets $S$ properly.
If we use this complete intersection to link $X$ to $W$, then by
proposition 1.11 we have an exact sequence
$$0 \ra \N^\vee(-s-t) \ra \P^\vee(-s-t) \oplus \O(-s) \oplus \O(-t) \ra \I_W
\ra 0$$
where the composition of the first map and projection to $\O(-t)$ is
precisely the map $\N^\vee \ra \O(p)$ tensored by $\O(-s-t)$, which is
surjective because the section is nowhere vanishing. Now we can form
a surjective morphism of exact sequences
$$\CD
0 @>>> \N^\vee(-s-t) @>>> \P^\vee(-s-t) \oplus\O(-s) \oplus \O(-t)@>>>
\I_W @. \rightarrow 0 \\
@. @VVV @VVV @VVV @. \\
0 @>>> \O(-t) @>>> \O(-t) @>>> 0 @.
\endCD $$
Applying the snake lemma gives an exact sequence of kernels
$$0 \ra {\N^\prime}^\vee(-s-t) \ra \P^\vee(-s-t) \oplus \O(-s) \ra \I_W \ra
0$$
We can choose a pair of equations in $I_W$ which are images of
generators of two of the line bundle summands in $\P^\vee(-s-t) \oplus \O(-s)$
and give hypersurfaces $S^\prime,T^\prime$ of degrees $s^\prime,t^\prime$
which meet properly.
If we apply proposition 1.8 we get a resolution for the
subscheme $X^\prime$ linked to $W$ by $S^\prime \cap T^\prime$.
Because of our choices of hypersurface, we can split off the
two line bundle summands corresponding to the equations of $S^\prime$
and $T^\prime$, which leaves an $\N$-type resolution for $X^\prime$ with
just $\N^\prime(s+t-s^\prime-t^\prime)$ in the middle.
\enddemo
\remark{Remark 2.6} Comparing with \cite{14, proposition 1.11}, we see
that the statement here is somewhat weaker, as the section of $\N$
considered is split instead of merely nowhere vanishing. The reason
is that the case $f=0$ in Rao's proof is flawed, but weakening
the hypotheses allows the proof to be repaired. Fortunately, the
weaker statement is still enough to prove proposition 2.8,
so Rao's theorem is not compromised.
\endremark
\proclaim{Lemma 2.7} Let $\N$ be a reflexive sheaf such that
$H^1_*(\N^\vee)=0$ and $\ext^1(\N,\O)=0$. Let $U \subset \Pn$
be an open subscheme whose complement has codimension $\geq 2$.
If $Y \subset U$ is a closed subscheme of codimension two (or empty)
which has a resolution of the form
$$0 \ra \P_U \ra \N_U \ra \I_{Y,U}(t) \ra 0$$
with $\P$ dissoci\'e, then there exists a subscheme ${\overline Y} \subset \Pn$
of pure codimension two which agrees with $Y$ on $U$ and has a resolution
$$0 \ra \P \ra \N \ra \I_{\overline Y}(t) \ra 0.$$
In particular, if the codimension of the complement of $U$ is $\geq 3$,
then ${\overline Y}$ is the scheme-theoretic closure of $Y$ in $\Pn$.
\endproclaim
\demo{Proof} Let $j:U \ra \Pn$ be the inclusion. Since both $\P$ and $\N$
are reflexive, we see from \cite{6, proposition 2.5} that $j_*(\P_U)=\P$
and $j_*(\N_U)=\N$. The left exactness of $j_*$ gives an exact sequence
$$0 \ra \P \ra \N \ra j_*(\I_{Y,U}(t)).$$
The left exactness of $j_*$ shows that $j_*(\I_{Y,U}(t))$ is a subsheaf
of $\O(t)$, so we get an exact sequence
$$0 \ra \P \ra \N \ra \I_{\overline Y}(t) \ra 0$$
where ${\overline Y}$ is a closed subscheme which agrees with $Y$ on
$U$. By corollary 1.18, this exact sequence shows that
${\overline Y}$ is of pure codimension two. It follows that $Y$ is as
well. If the complement of $U$ has codimension $\geq 3$, then ${\overline Y}$
can only be the scheme-theoretic closure of $Y$ (because $U$ must
contain all the generic points of each).
\enddemo
\proclaim{Proposition 2.8}
Let $X_1,X_2 \subset \Pn$ be two subschemes of pure
codimension two having $\N$-type resolutions of the form
$$0 \ra \oplus_{j=1}^r\O(a_j) @>\alpha>> \F(a) \ra \I_{X_1}
 \ra 0$$
$$0 \ra \oplus_{j=1}^r\O(b_j) @>\beta>> \F(b) \ra \I_{X_2}
 \ra 0$$
with $H^1_*(\F^\vee)=0$. Then $X_1,X_2$ belong to the same even linkage class.
\endproclaim
\demo{Proof} If $\F$ is dissoci\'e, then both $X_1$ and $X_2$ are
Arithmetically Cohen-Macaulay and we are done, so we may assume that $\F$ is
not dissoci\'e. In this case, proposition 2.3 shows that we can
split off line bundles from $\F$ until we obtain a representative $\F_0$ for
the stable equivalence class for $\F$ of minimal degree. Applying lemma
2.5 repeatedly allows us to replace $X_1,X_2$ with representatives
from their respective even linkage classes which have resolutions as above
involving $\F_0$. Thus we reduce to the case where $\F$ has no line bundle
summands. \par
Now the proof is identical to Rao's proof (with one slight change)
of the corresponding lemma for locally free sheaves \cite{15, lemma 1.6},
which we repeat for the reader's convenience. The maps
$\alpha$ and $\beta$ are described by two sequences $\{s_i\}$ and $\{t_i\}$ of
global sections of various twists of $\F$. We will assume that
$s_i=t_i$ for $i \leq j < r$ and show how to replace $X_1,X_2$ with
subschemes with resolutions in which we have $s_{j+1}=t_{j+1}$. \par
For this, we can find general sections $r_1,r_2$ of $\F(p)$ for $p >> 0$
whose images in $I_{X_1}$ give rise to independent hypersurfaces
$S_{11}$ and $S_{12}$ containing $X_1$ and whose images in $I_{X_2}$
give rise to independent hypersurfaces $S_{21}$ and $S_{22}$ containing
$X_2$. Using these hypersurfaces, we link $X_1$ to $Y_1$ and $X_2$ to
$Y_2$. By proposition 1.11, $Y_1$ has an $\E$-type resolution
$$0 \ra \F^\vee(a-2p) \ra \oplus_{i=1}^r\O(2a-2p-a_i) \oplus \O(a-p)^2 \ra
\I_{Y_2}
 \ra 0$$
where the restriction of the first map to the summand $\oplus\O(2a-2p-a_i)$ is
the dual map to $\alpha$. By our initial reduction, none of the sections
$s_i$ of $\F$ can be split off. Consider the induced map
$\O(2a-2p-a_{j+1}) \ra \I_{Y_1}$. If this map is zero, then the
map induced map $H^0(\O) \ra H^0(\I_{Y_1}(-2a+2p+a_{j+1}))$ is clearly
zero. Twisting the sequence above by $-2a+2p+a_{j+1}$, we see that
the map $H^0(\F^\vee(-a+a_{j+1})) \ra H^0(\O)$ is surjective, and we deduce
that the induced map $\F^\vee(-a+a_{j+1}) \ra \O$ is a split surjection.
It follows that the section $\O(-a_{j+1}) \ra \F(a)$ is also split,
which contradicts our original reduction. \par
It follows from the above paragraph that the map
$\O(2a-2p-a_{j+1}) \ra \I_{Y_1}$ is nonzero and hence defines a hypersurface
$T_1$ containing $Y_1$. Similarly the image of $\O(2b-2p-b_{j+1})$ in
$\I_{Y_2}$ gives a hypersurface $T_2$ which contains $Y_2$. If we replace
the sections $r_1,r_2$ with a general linear combination $lr_1+mr_2,hr_1+kr_2$,
we will get new pairs of hypersurfaces containing $X_1,X_2$ which form the
same complete intersections $S_{11} \cap S_{12}$ and $S_{21} \cap S_{22}$.
If we replace $r_1,r_2$ with a sufficiently general linear combination,
the corresponding hypersurface $S_{11}$ (respectively $S_{21}$) will meet $T_1$
(respectively $T_2$) properly. We use the corresponding complete intersection
$S_{11} \cap T_1$ (respectively $S_{21} \cap T_2$) to link $Y_1$
(respectively $Y_2$) to $X_1^1$ (respectively $X_2^1$). \par
By corollary 1.10, we obtain a resolution
$$0 \ra \oplus_{i \neq j+1}\O(a_i-a_{j+1}+a-p) \oplus \O(2a-2p-a_{j+1})
@>\alpha^1>> \F(2a-p-a_{j+1}) \ra \I_{X_1^1} \ra 0$$
where the map $\alpha^1$ is given by the sequence
$$\{s_1, \dots, s_j, s_{j+2}, \dots, s_r, -r_1\}$$
and we have effectively replaced $s_{j+1}$ with $r_1$. The analogous
statement holds for the resolution of $X_2^1$. We repeat this process
until we obtain subschemes $X_1^*,X_2^*$ in the respective even linkage
classes of $X_1,X_2$ such that the there are $\N$-type resolutions
with maps $\alpha^*,\beta^*$ which give the same sections of $\F$.
On the open set $U$ where $\F$ is locally free, we see that
$X_1^*=X_2^*$ because the ideals of both are given by the maximal
minors of $\alpha^*=\beta^*$. Now apply lemma 2.7 to see
that $X_1^*=X_2^*$. This shows that $X_1$ is evenly linked to $X_2$,
proving the proposition.
\enddemo
\proclaim{Corollary 2.9} The map $\Phi$ from proposition 2.4 is injective.
\endproclaim
\demo{Proof} Suppose that two schemes $Y_1,Y_2$ have $\E$-type
resolutions
$$0 \ra \E_i \ra \Q_i \ra \I_{Y_i} \ra 0$$
where $\E_1$ is stably equivalent to $\E_2$. By adding some extra
summands we may assume that $\E_1 \cong \E_2(d)$ for some $d \in \Z$.
If we choose two equations in each of the ideals $I_{Y_i}$ which come from
generators for line bundle summands of $\Q_i$ and define hypersurfaces
which meet properly, then we can link $Y_1,Y_2$ via these hypersurfaces
to $X_1,X_2$. Using proposition 1.8 we get $\N$-type
resolutions for the $X_i$, and by our choice of hypersurfaces, we
can split off the pair of summands corresponding to the hypersurfaces
used for the linkage. In this case, we get $\N$-type resolutions as
in the lemma, and hence we can conclude that $Y_1,Y_2$ are in the
same even linkage class.
\enddemo
\proclaim{Proposition 2.10}
The map $\Phi$ from proposition 2.4 is surjective.
\endproclaim
\demo{Proof} Once again we follow Rao's proof. Let $\E$ be a reflexive
sheaf such that $H^1_*(\E)=0$ and $\ext^1(\E^\vee,\O)=0$, and let
$\N=\E^\vee$. The set where $\N$ fails to be locally free is closed
of codimension $\geq 3$. Let $U$ be the open complement of this set
with inclusion map $j:U \hookrightarrow \Pn$. \par
If we choose $p$ large enough that $\N(p)$ is generated by global sections,
then $\N(p)_U$ is still generated by the same global sections. By
quotienting out by nowhere vanishing sections, we obtain a quotient
bundle $\M$ of rank $r \leq n$ which is still generated by global
sections. Now we can apply \cite{8, theorem 3.3} to see that a general
set of $r-1$ sections $\{s_i\}_{i=1}^{r-1}$ of $\M$ will have a
degeneracy locus $Y$ of pure codimension two (on $U$). The generators
of $\I_Y$ are locally given by the image of $\Lambda^{r-1}\gamma$,
where $\gamma:\O^{r-1} \ra \M$ is the map determined by the sections
$s_i$. It follows that $Y$ is locally Cohen-Macaulay with resolution
$$0 \ra \O_U^{r-1} \ra \M \ra \I_{Y,U}(t) \ra 0$$
for some twist $t$. Since the kernel of the map $\N(p)_U \ra \M$ was
$\O_U^l$, where $l$ is the difference in rank between $\N$ and $\M$,
we obtain an exact sequence
$$0 \ra \O_U^{r-1+l} \ra \N(p)_U \ra \I_{Y,U} \ra 0.$$
Now we apply lemma 2.7 to obtain an exact sequence
$$0 \ra \O^{r-1+l} \ra \N \ra \I_{\overline Y} \ra 0$$
where ${\overline Y}$ is of pure codimension two and in fact
the scheme-theoretic closure of $Y$. If we link ${\overline Y}$ to
another subscheme $X$, proposition 1.11
shows that $X$ has an $\E$-type resolution with last
step a twist of $\E$. This shows that $\Phi$ is surjective.
\enddemo
\proclaim{Theorem 2.11}
The even linkage classes of codimension two subschemes in $\Pn$ are
in bijective correspondence with the stable equivalence classes of
reflexive sheaves $\E$ on $\Pn$ satisfying $H^1_*(\E)=0$ and
$\ext^1(\E^\vee,\O)=0$.
\endproclaim
\proclaim{Corollary 2.12}
The even linkage classes of codimension two
subschemes in $\Pn$ are in bijective correspondence with the stable
equivalence classes of reflexive sheaves $\N$ on $\Pn$ satisfying
$H^1_*(\N^\vee)=0$ and $\ext^1(\N,\O)=0$ by taking the middle sheaf of
an $\N$-type resolution.
\endproclaim
\demo{Proof} This is evident from applying the mapping cone procedure
to $\E$-type resolutions in the dual class.
\enddemo
\remark{Remark 2.13}
In view of proposition 2.3, the stable equivalence classes of
reflexive sheaves are classified (up to twist) by certain reflexive
sheaves which have no line bundle summands.
If we think of stable equivalence classes as being determined by a
sheaf of minimal rank, we see that via $\E$-type resolutions, the
even linkage classes of codimension two subschemes are in bijective
correspondence with the set of reflexive sheaves $\E_0$ which have
no line bundle summands and satisfy $H^1_*(\E_0)$ and $\ext^1(\E^\vee,\O)=0$,
up to twist. In this case we will say that the even linkage class
$\L$ corresponds to the stable equivalence class $[\E_0]$ via
$\E$-type resolutions. Similarly, the even linkage classes of codimension two
subschemes are in bijective correspondence with the set of reflexive
sheaves $\N_0$ which have no line bundle summands and satisfy
$H^1_*(\N^\vee)=0$ and $\ext^1(\N,\O)=0$, up to twist. In this case
we say that the even linkage class $\L$ corresponds to the stable
equivalence class $[\N_0]$ via $\N$-type resolutions.
\endremark
\remark{Remark 2.14}
A little checking shows that proposition 2.4
generalizes to codimension $r$. In other words, by using the $\E$ in
the $\E$-type resolution, we get a well-defined map from even linkage
classes to stable equivalence classes of $rth$ syzygy modules. This map
is not a bijection in the case when $r=n$, so Rao's correspondence does
not extend to higher codimension without adaption.
\endremark
%
\heading The Lazarsfeld-Rao Property \endheading
\definition{Definition 3.1}
Let $X \subset \Pn$ be a subscheme of pure codimension two.
Let $S$ be a hypersurface of degree $s$ which contains $X$ and let $h \geq 0$
be an integer. We say that $Y$ is obtained from $X$ by a {\it basic double
link of height} $h$ {\it on} $S$ if there is a hypersurface $H$ of degree
$h$, a hypersurface $T$ and a subscheme $Z$ such that $X$ is linked to $Z$
by $S \cap T$ and $Z$ is linked to $Y$ by $S \cap (T \cup H)$.
\enddefinition
\definition{Definition 3.2}
Let $X,X^\prime$ be subschemes of pure codimension two
in $\Pn$. We say that $X^\prime$ {\it dominates} $X$ {\it at height}
$h \geq 0$ if $X^\prime$ can be obtained from $X$ by a sequence of basic
double links with heights summing to $h$, followed by a deformation which
preserves cohomology and even linkage class. In this case we write
$X \leq_h X^\prime$, or simply $X \leq X^\prime$ if $h$ is not specified.
\enddefinition
We will show that domination is transitive and hence defines a partial
ordering on an even linkage class. Moreover, we will characterize
domination in terms of the $\N$-type resolutions discussed in section 2.
\proclaim{Lemma 3.3}
Let $\E @>\phi>> \Q$ be an injective morphism of reflexive sheaves
of the same degree on $\Pn$, where $\E$ has rank $r$ and $\Q$ has rank $r+1$.
Assume that $\E$ and $\Q$ satisfy one of the following sets of conditions:
\roster
\item $H^1_*(\E)=0$, $\ext^1(\E^\vee,\O)=0$ and $\Q$ is dissoci\'e.
\item $H^1_*(\Q^\vee)=0$, $\ext^1(\Q,\O)=0$ and $\E$ is dissoci\'e.
\endroster
Then the following are equivalent. \roster
\item The cokernel of $\phi$ is isomorphic to the ideal sheaf of a subscheme
of pure codimension two.
\item The cokernel of $\phi$ is torsion free.
\item The complex $$0 \ra \E @>\phi>> \Q @>\Lambda^r\phi^\vee \otimes 1>> \O$$
is exact. \endroster
\endproclaim
\demo{Proof} First a word of explanation about the complex from condition (3).
We have an isomorphism $\Q \cong (\Lambda^r\Q)^\vee \otimes \Lambda^{r+1}\Q$
because these are isomorphic on the open set $U$ whose complement has
codimension $\geq 2$ and both sides are reflexive. Since both $\Q$ and $\E$
have the same degree, $(\Lambda^r\E)^\vee \otimes \Lambda^{r+1}\Q \cong \O$
which shows explains the map $(\Lambda^r\phi)^\vee \otimes 1$. To see that
the sequence of condition (3) is a complex, first note that this holds on
the open set $U$. If $j:U \hookrightarrow \Pn$ is the inclusion map, then
applying $j_*$ to the sequence on $U$ gives the sequence from condition (3)
because all three sheaves involved are reflexive. Since $j_*$ is a functor,
the resulting sequence is a complex. \par
Now we may write the image of $(\Lambda^r\phi)^\vee \otimes 1$ as
$\I(-\alpha)$, where $\I$ is the ideal sheaf of a subscheme in $\Pn$ of
codimension $\geq 2$ and $\alpha \geq 0$. Letting $\E^\prime$ be the
kernel of the induced map $\Q \ra \I(-\alpha)$, we have a commutative
diagram
$$\CD
0 @>>> \E @>>> \Q @>>> cok\phi @>>> 0 \\
@. @VVV @VVV @VVV @. \\
0 @>>> \E^\prime @>>> \Q @>>> \I(-\alpha) @>>> 0
\endCD $$
where the rows are short exact sequences. Clearly the cokernel of the
inclusion $\E^\prime \ra \E$ is a torsion sheaf, which is isomorphic
to the kernel of the surjection $cok\phi \ra \I(-\alpha)$ by the snake
lemma. Now we show the equivalence of the conditions. \par
$(1) \Rightarrow (2)$ This is trivial. \par
$(2) \Rightarrow (3)$ Since the cokernel of $\phi$ is torsion free, the torsion
sheaf $\E/{\E^\prime}$ must be zero, when the complex of (3) is exact. \par
$(3) \Rightarrow (1)$ If the complex is exact, then the cokernel of $\phi$ is
an ideal sheaf. The associated short exact sequence shows that the degree of
this ideal sheaf is zero, hence the subscheme $Z$ that it defines is of
codimension $\geq 2$. Now in either of the cases (1) or (2) given above,
we can apply either corollary 1.18 or corollary 1.16 to
see that $Z$ is of pure codimension two. \par
\enddemo
\proclaim{Lemma 3.4}
Let $X @>f>> Y$ be a morphism of schemes, and
let $0 \ra A \ra B \ra C \ra 0$ be a complex of sheaves on $X$. Assume
further that both $B$ and $C$ are flat over $Y$. If $y=f(x)$, then the set
$$\{x \in X:0 \ra A_x \otimes k(y) \ra B_x \otimes k(y) \ra C_x \otimes k(y)
\text{is exact}\}$$
is open in $X$.
\endproclaim
\demo{Proof} Let $I$ and $J$ denote the image and cokernel of the map $A \ra
B$.
Since the sequence of the lemma is a complex, it can be factored as in the
following commutative diagram:
$$\CD
@. A @. @. @. \\
@. @VVV @. @. @. \\
0 @>>> I @>>> B @>>> J @>>> 0 \\
@. @. @. @VVV @. \\
@. @. @. C @.
\endCD $$
This diagram remains exact when we localize at $x$. Applying the right
exact functor $\otimes k(y)$, we see that the sequence
$$0 \ra A_x \otimes k(y) \ra B_x \otimes k(y) \ra C_x \otimes k(y)$$
is exact if and only if the maps $A_x \otimes k(y) \ra B_x \otimes k(y)$
and $J_x \otimes k(y) \ra C_x \otimes k(y)$ are injective. Since $B$ and $C$
are both flat over $Y$, the set of $x \in X$ such that both these maps are
injective is open by \cite{3, IV, corollary 11.3}.
\enddemo
\proclaim{Corollary 3.5}
Let $X @>f>> Y$ be a proper flat morphism of
noetherian schemes, and let $0 \ra A \ra B \ra C \ra 0$ be a complex
on $X$ with $B$ and $C$ locally free. Then the set
$$U=\{y \in Y: 0 \ra A \otimes k(y) \ra B \otimes k(y) \ra C \otimes k(y)
\text{is exact}\}$$
is open in $Y$.
\endproclaim
\demo{Proof} Let $y \in Y$. Then $y \in U$ if and only if the sequence
$0 \ra A_x \otimes k(y) \ra B_x \otimes k(y) \ra C_x \otimes k(y)$ is
exact for each $x \in f^{-1}(y)$. This condition is open on $X$ by the
previous lemma, because both $B$ and $C$ are flat over $Y$ in this case.
If $y \in U$, then we can find an open set $V_x$ for each $x \in f^{-1}(y)$
giving the exactness condition. Since $f^{-1}(y)$ is noetherian, a finite
number of these open sets cover $f^{-1}(y)$, call them $\{V_i\}_1^n$.
Since $f$ is proper, the set $f(X-\cup{V_i})$ is closed in $Y$ and avoids $y$.
Hence there is an open neighborhood $W$ of $y$ which also avoids the set
$f(X-\cup{V_i})$. Since for any $w \in W$ we have that
$f^{-1}(w) \subset \cup{V_i}$, we see that $W \subset U$. This shows that
$U$ is open.
\enddemo
\proclaim{Lemma 3.6}
Let $\L$ be an even linkage class of subschemes
of pure codimension two in $\Pn$.
Let $X,Y \in \L$ such that $h^i(\I_X(n))=h^i(\I_Y(n))$ for all
$n \in \Z$ and $i \geq 0$. Then there exists an irreducible flat family of
subschemes $\{X_s\}_{s \in S}$ with constant cohomology from $\L$
to which $X$ and $Y$ belong.
\endproclaim
\demo{Proof} Since $h^0(\I_X(n))=h^0(\I_Y(n))$ for all $n$, we can find a
graded free $S$-module $Q$ with surjections $Q \ra I_X$ and $Q \ra I_Y$.
Sheafifying gives $\E$-type resolutions
$$0 \ra \E_X \ra \Q \ra \I_X \ra 0$$
$$0 \ra \E_Y \ra \Q \ra \I_Y \ra 0$$
where the kernels $\E_X$ and $\E_Y$ are stably equivalent by theorem 2.11.
If we write $\E_X \oplus \P_1 \cong \E_Y(d) \oplus \P_2$, we consider two
cases. If all the middle cohomology of both of these sheaves is zero, then
$\E_X$ and $\E_Y$ are both dissoci\'e. In this case, the fact that
$h^0(\I_X(n))=h^0(\I_Y(n))$ for all $n$ shows that $\E_X \cong \E_Y$.
If some of the middle cohomology is nonzero, then the fact that
$h^i(\I_X(n))=h^i(\I_Y(n))$ for $1 \leq i < n-1$ shows that $d=0$.
In this case we can add $\P_1$ to the first sequence and $\P_2$ to the
second sequence to get new $\E$-type resolutions in which the leftmost
sheaves are isomorphic and the middle dissoci\'e sheaves are
$\Q \oplus \P_1$ and $\Q \oplus \P_2$. Since the cohomology of the middle
sheaves must be the same and they are both dissoci\'e, we see that they
must be isomorphic. In either case, we obtain $\E$-type resolutions
$$0 \ra \E \ra \Q \ra \O$$
for $X$ and $Y$ which involve the same sheaves (but possibly different maps).
\par
Now let $U$ be the subset of $\Pn$ where $\E$ is locally free.
By \cite{6, corollary 1.4}, the codimension of the complement of $U$
is $\geq 3$. Let $V$ denote the parameter space for all homomorphisms
$\phi:\E \ra \Q$. Then if $p:\Pn \times V \ra \Pn$ is the natural projection,
there is a universal homomorphism $\phi:p^*\E \ra p^*\Q$. On the open
subset $U \times V \subset \Pn \times V$ we have the complex
$$0 \ra p^*\E @>\phi>> p^*\Q @>\Lambda^r\phi^\vee \otimes 1>> p^*\O \ra 0$$
If $j$ is the inclusion map for this open subset, applying $j_*$
extends the complex to all of $\Pn \times V$. By the preceding corollary,
the set $V_1$ of $v \in V$ such that
$0 \ra p^*\E \otimes k(v) \ra p^*\Q \otimes k(v) \ra p^*\O \otimes k(v)$
is exact is open. Moreover, the morphisms corresponding to $X$ and $Y$ show
that $V_1$ is nonempty. Since $V$ is smooth and irreducible, so is the nonempty
open subset $V_1$. $V_1$ gives a family of subschemes $X_v$ with exact
sequences
$$0 \ra \E @>\phi_v>> \Q \ra \O \ra \O_{X_v} \ra 0$$
which specialize to $X$ and $Y$. These subschemes are of pure codimension
two by lemma 3.3. Theorem 2.11 shows that they all lie in $\L$,
which proves the statement of the lemma.
\enddemo
\definition{Definition 3.7}
If $f:\Z \ra {\Bbb N}$ is a function such that $f(n)=0$ for
$n << 0$, then we define the function $f^\#$ by $f^\#(a)=\sum_{n \leq a}f(n)$.
\enddefinition
\remark{Remark 3.8}
We will be interested in comparing two such
functions which have finite support and satisfy $\sum{f(n)}=\sum{g(n)}$.
In this case, let $\{a_i\}$ and $\{b_i\}$ be the sequences of integers
defined by the conditions $f^\#(a_i-1) < i \leq f^\#(a_i)$ and
$g^\#(b_i-1) < i \leq g^\#(b_i)$. These are the unique nondecreasing sequences
that satisfy $f(n)=\#\{i:a_i=n\}$ and $g(n)=\#\{i:b_i=n\}$. Note that
$\#\{a_i\}=\sum{f(n)}=\sum{g(n)}=\#\{b_i\}$. Now an elementary calculation
shows that $f^\#(n) \geq g^\#(n) \iff a_i \leq b_i$ for each $i$. This
characterization will be used in the next proposition.
\endremark
\proclaim{Proposition 3.9}
Let $X,Y \subset \Pn$ be subschemes of pure
codimension two having $\N$-type resolutions
$$0 \ra \oplus\O(-n)^{r(n)} \ra \N \ra \I_X(h_1) \ra 0$$
$$0 \ra \oplus\O(-n)^{s(n)} \ra \N \ra \I_Y(h_2) \ra 0$$
such that $H^1_*(\N^\vee)=0$. Then $X \leq Y \iff r^\# \geq s^\#$.
\endproclaim
\demo{Proof} Suppose that $X \leq Y$. First we show that there {\it exist}
$\N$-type resolutions as above. If $Y$ is obtained from $X$ by a basic
double link of type $(s,h)$ and $X$ has an $\N$-type resolution
$$0 \ra \P \ra \N \ra \I_X \ra 0$$
then $Y$ has an $\N$-type resolution of the form
$$0 \ra \P \oplus \O(-s) \ra \N \oplus \O(-s+h) \ra \I_Y(h) \ra 0.$$
If we add the summand $\O(-s+h)$ to both terms in the resolution for
$\I_X$, we see that the dissoci\'e part becomes $\P \oplus \O(-s+h)$.
Comparing this to $\P \oplus \O(-s)$, we see that $r^\# \geq s^\#$ for
this case.  The same argument shows that if we do
a sequence of basic double links to get from $X$ to $Y$, then we can
find such resolutions of $X$ and $Y$. \par
Now we assume that $X \leq Y$, so that $Y$ can be obtained from
$X$ by a sequence of basic double links followed by a deformation which
preserves cohomology and even linkage class. Let $Z$ be the subscheme obtained
from $X$ by the sequence of basic double links, so that we can write
sequences
$$0 \ra \oplus\O(-n)^{r(n)} \ra \N \ra \I_X \ra 0$$
$$0 \ra \oplus\O(-n)^{s(n)} \ra \N \ra \I_Z(h) \ra 0$$
with $r^\# \geq s^\#$. Let $Y$ be obtained from $Z$ by a deformation
which preserves cohomology and even linkage class. Since $Y$ is in
the same even linkage class as $X$, we can obtain $Y$ from $X$ by
an even sequence of links. Applying the mapping cone procedure,
we get a resolution for $Y$ of the form
$$0 \ra \P \ra \N \oplus \Q \ra \I_Y(h) \ra 0$$
where $\P,\Q$ are dissoci\'e (if some of the middle cohomology of $\N$ is
nonzero, it's clear that the twist $h$ is the same for both $Z$ and $Y$.
If all the middle cohomology vanishes, then $\N$ is dissoci\'e and it's
also clear that we can use the same twist $h$ for both $Z$ and $Y$).
Add the dissoci\'e sheaf $\Q$ to the $\N$-type resolutions for $X$ and $Z$
above. With these new resolutions, we still have the
inequalities $r_1^\# \geq s_1^\#$,
where $\oplus\O(-n)^{s_1(n)} \cong \oplus\O(-n)^{s(n)} \oplus \Q$ and
$\oplus\O(-n)^{r_1(n)} \cong \oplus\O(-n)^{r(n)} \oplus \Q$. Further,
the fact that the cohomology of $Y$ and the cohomology of $Z$ are
the same shows that $\P \cong \oplus\O(-n)^{s_1(n)}$. This proves
that there exist sequences as in the statement of the proposition. \par
To see that the inequalities must hold for {\it any} such $\N$-type
resolution, we show that they are independent of the $\N$-type resolution
chosen. Since the $\N$'s are stably equivalent, it suffices to show this
assuming that the $\N$'s differ by a single line bundle $\O(-a)$. This is clear
because in adding this line bundle the corresponding functions
$r^\#(n)$ and $s^\#(n)$ are both increased by $1$ for all $n \geq a$ to get
the new functions $r_1^\#$ and $s_1^\#$. It follows that $r^\# \geq s^\# \iff
r_1^\# \geq s_1^\#$. \par
   Now assume that $X$ and $Y$ have such $\N$-type resolutions.
We wish to show that $X \leq Y$. This is essentially \cite{10, IV,
lemma 5.2}. We give a slightly different proof here. Rewrite
$\oplus\O(-n)^{r(n)}=\oplus\O(-a_i)$ and $\oplus\O(-n)^{s(n)}=\oplus\O(-b_i)$
with the $a_i$ and $b_i$ nondecreasing. By the preceding remark, we have
that $a_i \leq b_i \iff r^\# \geq s^\#$. Let $h=\sum{b_i}-\sum{a_i}$
By looking at Chern classes, we can see that $h=h_2-h_1$. I'll first show by
induction on $h$ that there exists a scheme $Y^\prime$ which is obtained
from $X$ by a sequence of basic double links which has the same
cohomology and Rao modules as $Y$. For the induction base case
$h=0$ we have $a_i=b_i$ for all $i$ and hence $r=s$, when the exact
sequences show that $h^0(\I_X(n))=h^0(\I_Y(n))$ for all $n$ and that
the higher Rao modules of $X$ and $Y$ agree, so taking a basic double
link of height $0$ will do (in other words, take $Y^\prime=X$). \par
Now consider the case where $h>0$. Let $n=max\{i:b_i>a_i\}$ so that
$a_i=b_i$ for $i>n$.
We claim that $h^0(\I_X(h_1+n)) \neq 0$. From the exact sequence for
$\I_X$, we see that
$$h^0(\I_X(h_1+n))=h^0(\N(n))-\sum_{a_i \leq n}{h^0(\O(n-a_i))}$$
but this is strictly greater than
$$h^0(\N(n))-\sum_{b_i \leq n}{h^0(\O(n-b_i))}=h^0(\I_Y(h_2+n)) \geq 0$$
and this proves the claim. This allows us to obtain a new scheme
$X_1$ from $X$ by a basic double link of type $(h_1+n,b_n-a_n)$.
$X_1$ has a resolution of the form
$$0 \ra \oplus\O(-a_i) \oplus \O(-b_n) \ra \N \oplus \O(-a_n) \ra
 \I_{X_1}(h_1+b_n-a_n) \ra 0$$
If we add an extra line bundle summand $\O(-a_n)$ to the resolution for
$Y$ we get
$$0 \ra \oplus\O(-b_i) \oplus \O(-a_n) \ra \N \oplus \O(-a_n) \ra
 \I_Y(h_2) \ra 0$$
If we reorder the summands of these, we find that the inequalities
still hold. If we sum the differences, we find
that the new value of $h$ is $h-b_n+a_n<h$, so we can apply the induction
hypothesis to obtain a subscheme $Y^\prime$ which is obtained by a
sequence of basic double links from $X_1$ and has the same cohomology
and higher Rao modules as $Y$. It follows that $Y^\prime$ can be obtained
from $X$ by a sequence of basic double links. Now we can apply lemma
3.6 to see that $Y^\prime$ can be deformed with constant
cohomology and even linkage class to $Y$. Thus $X \leq Y$.
\enddemo
\proclaim{Corollary 3.10}
Domination is transitive. In particular, it defines a
partial ordering on an even linkage class $\L$.
\endproclaim
\demo{Proof} If $X \leq Y$ and $Y \leq Z$, then we can find $\N$-type
resolutions as in proposition 3.9.
Using the stable equivalence of the  reflexive sheaves in the middle,
we can add line bundle summands to the three resolutions to obtain
exact sequences
$$0 \ra \oplus\O(-n)^{r(n)} \ra \N \ra \I_X \ra 0$$
$$0 \ra \oplus\O(-n)^{s(n)} \ra \N \ra \I_Y(h_1) \ra 0$$
$$0 \ra \oplus\O(-n)^{t(n)} \ra \N \ra \I_Z(h_2) \ra 0$$
with $r^\# \geq s^\#$ and $s^\# \geq t^\#$. It follows that $r^\# \geq t^\#$,
so proposition 3.9 shows that $X \leq Z$.
\enddemo
\definition{Definition 3.11}
Let $\L$ be an even linkage class of subschemes of pure
codimension two in $\Pn$. We say that $\L$ has the
{\it Lazarsfeld-Rao property} if $\L$ has
a minimal element with respect to domination.
\enddefinition
\remark{Remark 3.12}
In view of the definition of domination and the fact that domination is
a partial order, this definition is equivalent to the usual definition
used in \cite{10} or \cite{2}.
\endremark
\remark{Remark 3.13}
   In view of proposition 3.9 we can try to
produce minimal subschemes for an even linkage class $\L$ by finding
injections $\oplus\O(-n)^{r(n)} \hookrightarrow \N$ whose cokernel is a twisted
ideal sheaf of a subscheme of pure codimension two ($\N$ is in the stable
equivalence class corresponding to the $\N$-type resolutions) such that
the function $r^\#$ is as large as possible. Specifically, it would be ideal to
have a function $q$ associated to the reflexive sheaf $\N$ such that there
exists an injection $\oplus\O(-n)^{q(n)} \ra \N$ as above and has the property
that if $\oplus\O(-n)^{r(n)} \ra \N$ is another such injection, then
$r^\#(n) \leq q^\#(n)$ for each $n \in \Z$. Martin-Deschamps and Perrin have
already carried out this work in the case when $\N$ is a vector bundle
on $\Pthree$. Fortunately, their results about injections
$\oplus\O(-n)^{r(n)} \hookrightarrow \N$ are stated for the case when $\N$
is only a torsion free sheaf on $\Pthree$ and their proofs work on $\Pn$.
In particular, their solution applies to the reflexive sheaves $\N$
that come from $\N$-type resolutions. What follows is a brief summary of
their results, stated in the generality in which their proofs hold.
\endremark
\definition{Definition 3.14}
Let $B \subset A$ be an inclusion of $\O$-modules.
We say that $B$ is {\it maximal in} $A$ if for each $\O$-module $B^\prime$
having the same rank as $B$ and satisfying $B \subset B^\prime \subset A$
we have $B = B^\prime$.
\enddefinition
\proclaim{Proposition 3.15}
Let $B \subset A$ be an inclusion of $\O$-modules,
where $B$ is locally free and $A$ is torsion free.
Then the following are equivalent: \roster
\item $B$ is maximal in $A$.
\item $A/B$ is torsion free.
\item $A/B$ is torsion free in codimension one.
\item $A/B$ is locally free in codimension one.
\item $A/B$ has constant rank in codimension one.
\item $A/B$ is locally a direct factor of $A$ in codimension one.
\endroster
\endproclaim
\demo{Proof} See \cite{10, IV, proposition 1.2}. In fact, their proof shows
that these are equivalent if $B$ is only reflexive.
\enddemo
\definition{Definition 3.16}
Let $\N$ be a torsion free sheaf on $\Pn$. Let
$\oplus{S(-n)^{l(n)}} \ra H^0_*(\N)$ be a minimal graded surjection of
$S$-modules. This defines the function $l_{\N}=l$. Sheafifying this
surjection defines the map $\sigma:\oplus\O(-n)^{l(n)} \ra \N$.
For $a \in \Z$, the map $\sigma_a$ is defined by restricting $\sigma$
to $\oplus_{n \leq a}\O(-n)^{l(n)}$.
\enddefinition
\remark{Remark 3.17} For each $a \in \Z$, $\sigma_a$ is well-defined up to a
choice of basis for
$\oplus_{n \leq a}\O(-n)^{l(n)}$ and hence the image of $\sigma_a$
(which is denoted by $\N_{\leq a}$ in \cite{10}) depends only on $a$.
The function $l=l_{\N}$ could have been defined by the formula
$l(n)=dim_k(N \otimes_S k)_n$, where $N=H^0_*(\N)$.
\endremark
\definition{Definition 3.18} Let $\N$ be a torsion free sheaf on $\Pn$ and let
$\sigma_a$ be as in definition 3.16. Then we define
$$ \alpha_a=rank(\sigma_a) $$
$$ \beta_a=inf_D(rank(\sigma_a|_D)) $$
where this infimum is taken over all integral divisors $D$. We always have
the inequalities
$$0 \leq \beta_a \leq \alpha_a \leq l^\#(a)$$
and we define $a_0=sup\{a \in \Z:\alpha_a=\beta_a=l^\#(a)\}$ if this exists,
otherwise $a_0=\infty$. We also define
$a_1=inf\{a \in \Z:\alpha_a=\beta_a=rank(\N)\}$.
\enddefinition
\definition{Definition 3.19}
Let $\N$ be a torsion free sheaf on $\Pn$. Then we define
the function $q=q_\N$ by \roster
\item $q(n)=l(n)$ if $n < a_0$
\item $q^\#(n)=inf\{\alpha_n-1,\beta_n\}$ otherwise.
\endroster
\enddefinition
\remark{Remark 3.20}
This is the function referred to in remark 3.13. It is a consequence
of the results which follow that if $\N$ is reflexive with $\ext^1(\N,\O)=0$
and $H^1_*(\N^\vee)=0$, then there is an injection $\oplus\O(-n)^{q(n)} \ra \N$
whose cokernel is the twisted ideal sheaf of a subscheme of pure codimension
two. Further, for any other such injection $\oplus\O(-n)^{r(n)} \ra \N$
we have that $r^\#(n) \leq q^\#(n)$ for all $n \in \Z$.
\endremark
\proclaim{Proposition 3.21}
Let $\N$ be a torsion free sheaf on $\Pn$ with
associated functions $l,q,a_0$ and $a_1$ as in the preceding definitions.
\roster
\item If $\N$ is dissoci\'e, then $a_0=\infty$, $q=l$ and
$q^\#(\infty)=rank(\N)$.
\item If $\N$ is not dissoci\'e, then $q^\#(a_1)=q^\#(\infty)=rank(\N)-1$,
$0 \leq q(n) \leq l(n)$ for each $n \in \Z$, $a_0 \leq a_1 < \infty$,
$q(a_0) < l(a_0)$ and $q(n)=0$ for $n>a_1$.
\endroster
\endproclaim
\demo{Proof} See \cite{10, IV, proposition 2.7}.
\enddemo
\proclaim{Proposition 3.22}
Let $\N$ be a torsion free sheaf on $\Pn$ and suppose that
$\N^\prime = \N \oplus \P$ where $\P$ is dissoci\'e. Then
$q_{\N^\prime}=q_\N+l_\P$.
\endproclaim
\demo{Proof} See \cite{10, IV, proposition 2.9}.
\enddemo
\proclaim{Theorem 3.23}
Let $\N$ be a torsion free sheaf on $\Pn$ with
associated function $q=q_\N$. If an injection
$\oplus\O(-n)^{r(n)} \hookrightarrow \N$ has image maximal in $\N$,
then $r^\# \leq q^\#$.
\endproclaim
\demo{Proof} See \cite{10, IV, theorem 3.1}.
\enddemo
\proclaim{Theorem 3.24}
Let $\N$ be a torsion free sheaf on $\Pn$ with
associated function $q=q_\N$ and $\sigma:\oplus\O(-n)^{l(n)} \ra \N$. Then
$\oplus\O(-n)^{l(n)}$ has a direct factor of the form $\oplus\O(-n)^{q(n)}$
whose image under $\sigma$ is maximal in $\N$. Conversely, if
$\oplus\O(-n)^{q(n)} @>\beta>> \N$ is injective with maximal image,
then there is a split injection
$\tau:\oplus\O(-n)^{q(n)} \ra \oplus\O(-n)^{l(n)}$ such that
$\beta=\sigma \circ \tau$.
\endproclaim
\demo{Proof} For the first part, see \cite{10, IV, theorem 3.4}. The
converse part is the second part of \cite{10, IV, theorem 3.7}.
\enddemo
\proclaim{Proposition 3.25}
Let $\N$ be a torsion free sheaf on $\Pn$ with
associated function $q=q_\N$ and let $\P$ be dissoci\'e. If
$\tau:\oplus\O(-n)^{q(n)} \oplus \P \ra \N \oplus \P$ is an injection
with maximal image, then there is an automorphism $\nu$ of
$\oplus\O(-n)^{q(n)} \oplus \P$ such that the diagram
$$\CD
\oplus\O(-n)^{q(n)} \oplus \P @>\tau \circ \nu>> \N \oplus \P \\
@VVV @VVV \\
\P @>>> \P
\endCD $$
commutes, where the vertical arrows are the canonical projections.
\endproclaim
\demo{Proof} See \cite{10, IV, proposition 3.9}.
\enddemo
\proclaim{Theorem 3.26}
Let $\L$ be an even linkage class corresponding to the
stable equivalence class $[\N_0]$ via $\N$-type resolutions.
Assume $\N_0 \neq 0$ and let $q=q_{\N_0}$ be the associated function.
Then there exists $X_0 \in \L$ having an $\N$-type resolution of the form
$$0 \ra \oplus\O(-n)^{q(n)} \ra \N_0 \ra \I_{X_0}(h) \ra 0$$
and satisfying $X_0 \leq Y$ for each $Y \in \L$. In particular, $\L$ has
the Lazarsfeld-Rao property. $X_0^\prime$ is another such minimal element
of $\L \iff X_0^\prime$ has a resolution of the same
form $\iff X_0$ deforms to $X_0^\prime$ with constant cohomology through
subschemes in $\L$.
\endproclaim
\demo{Proof} There exists an injective homomorphism
$\tau:\oplus\O(-n)^{q(n)} \ra \N_0$ whose image is maximal in $\N_0$ by
theorem 3.24. From proposition 3.21 we have that
$q^\#(a_1)=rank(\N_0)-1$, so
the cokernel of $\tau$ has rank one. By maximality, this cokernel is
torsion free, hence it is a twisted ideal sheaf of a subscheme in $\Pn$
of pure codimension two by lemma 3.3 (after we twist the morphism
until the degrees become equal, condition (2) of lemma 3.3 holds). Let this
subscheme be $X_0$. \par
Now let $Y \in \L$. Then $Y$ has an $\N$-type resolution of the form
$$0 \ra \oplus\O(-n)^{r(n)} @>\nu>> \N_0 \oplus \Q \ra \I_Y(k) \ra 0$$
with $\Q$ dissoci\'e.  Since $\I_Y(k)$ is torsion free, the image of $\nu$ is
maximal in $\N_0 \oplus \Q$. It follows from theorem 3.23 that if
$q_1$ is the function associated to $\N_0 \oplus \Q$,
then we have $r^\# \leq q_1^\#$.
On the other hand, if we add the trivial summand $\Q$ on to the resolution
for $X_0$ and apply proposition 3.22, the resolution for $X_0$ corresponds
to this function $q_1$. In view of proposition
3.9 and remark 3.20, we see that $X_0 \leq Y$. \par
Finally, let $X_0^\prime$ be another subscheme for $\L$.
If $X_0^\prime$ has a resolution of the form that $X_0$ does, then
$h^i(\I_{X_0^\prime}(n))=h^i(\I_{X_0}(n))$ for all $n,i$, and hence by
lemma 3.6 we see that $X_0$ can be deformed to $X_0^\prime$
with constant cohomology and even linkage class. If there exists
such a deformation, then by definition of $\leq$ we have that
$X_0^\prime \leq_0 X_0$, hence $X_0^\prime$ is also minimal for $\L$.
It remains to be shown that if $X_0^\prime$ is minimal, then it has a
resolution of the same form that $X_0$ does. \par
If we write an $\N$-type resolution for $X_0^\prime$ as
$$0 \ra \oplus\O(-n)^{r(n)} \ra \N_0 \oplus \Q \ra \I_{X_0^\prime}(l) \ra 0$$
and add the trivial summand $\Q$ to the resolution for $X_0$, we conclude
that $q_1^\# \leq r^\#$ by minimality of $X_0^\prime$ and that
$r^\# \leq q_1^\#$ from theorem 3.23. It follows that $q_1^\#=r^\#$.
Since these are both functions of finite support, we see that $q_1=r$.
Applying proposition 3.22, we have that $r=q+l_\Q$. Now we apply
proposition 3.25 to change the map so that the image in
$\N_0 \oplus \Q$ is the same, but in such a way that the summand $\Q$
can be split off. Splitting off $\Q$, we get a resolution
$$0 \ra \oplus\O(-n)^{q(n)} \ra \N_0 \ra \I_{X_0^\prime} \ra 0.$$
This completes the proof.
\enddemo
\remark{Remark 3.27}
If $X_0$ is minimal for $\L$ as in the above theorem,
then there exists an exact sequence
$$0 \ra \oplus\O(-n)^{q(n)} \ra \N \ra \I_{X_0}(h) \ra 0$$
Since $q(n)=l(n)$ for $n<a_0$ and $q(a_0)<l(a_0)$ by proposition
3.21, we have that $$s(X_0)=min\{n:h^0(\I_{X_0}(n)) \neq 0\}=a_0+h.$$
\endremark
\proclaim{Proposition 3.28}
Let $\L$ be the even linkage class corresponding
to the stable equivalence class $[\N_0]$ via $\N$-type resolution, with
$\N_0 \neq 0$. Let
$\sigma:\oplus\O(-n)^{l(n)} \ra \N_0$ be the surjection of definition
3.16 and let $\E_0$ denote the kernel of $\sigma$. Then $X_0$ is
minimal for $\L \iff$ $X_0$ has an $\E$-type resolution of the form
$$0 \ra \E_0 \ra \oplus\O(-n)^{l(n)-q(n)} \ra \I_{X_0}(h) \ra 0.$$
\endproclaim
\demo{Proof} Theorem 3.26 gives an $\N$-type resolution for the
minimal subscheme $X_0$. By the second part of theorem 3.24 the
corresponding map $\oplus\O(-n)^{q(n)} \hookrightarrow \N_0$ factors
through $\sigma$ as a direct factor of $\oplus\O(-n)^{l(n)}$. It follows
that we have a commutative diagram
$$\CD
0 @>>> \oplus\O(-n)^{q(n)} @>>> \oplus\O(-n)^{l(n)} @>>>
\oplus\O(-n)^{l(n)-q(n)} @>>> 0 \\
@. @VVV @VVV @VVV @. \\
0 @>>> \oplus\O(-n)^{q(n)} @>>> \N_0 @>>> \I_{X_0}(h) @>>> 0
\endCD $$
whose rows are exact sequences and the leftmost vertical map is an isomorphism.
Since the middle map is surjective with kernel $\E_0$, the same holds for the
rightmost map by the snake lemma, which gives the resolution. Because
the induced map $H^0_*(\oplus\O(-n)^{l(n)}) \ra H^0_*(\I_{X_0}(h))$ is
surjective, this also holds for the map $H^0_*(\oplus\O(-n)^{l(n)-q(n)}) \ra
H^0_*(\I_{X_0}(h))$ and hence the resolution is actually an $\E$-type
resolution. \par
Conversely, if $X_0^\prime$ has such a resolution, it is in the same
even linkage class as some minimal subscheme $X_0$ for $\L$ by
theorem 2.11 and this resolution shows that $X_0^\prime$ has the
same cohomology as $X_0$. Applying lemma 3.6 we see that
$X_0$ can be deformed with constant cohomology and even linkage class to
$X_0^\prime$, hence $X_0^\prime$ is also minimal by theorem 3.26.
\enddemo
\proclaim{Proposition 3.29}
Let $\E_0$ be the kernel of the map
$\sigma:\oplus\O(-n)^{l(n)} \ra \N_0$ as in proposition 3.28.
Let $q$ be the function associated to $\N_0$ and let $q^\prime$ be the
function associated to $\E_0^\vee$. Then $q$ and $q^\prime$ are related
by the formula
$$q(n)+q^\prime(n)=l(n)-\epsilon_0(n)-\epsilon_1(n)$$
where $\epsilon_0(a_0)=1$, $\epsilon_0(n)=0$ for $n \neq 0$,
$\epsilon_1(a_1)=1,$ and $\epsilon_1(n)=0$ for $n \neq a_1$.
\endproclaim
\demo{Proof} See \cite{10, IV, proposition 5.7}.
\enddemo
\proclaim{Lemma 3.30}
Let $f_i$ be homogeneous polynomials of nondecreasing degrees
$d_1 \leq d_2 \leq \dots \leq d_r$ which cut out a subscheme in $\Pn$
of codimension $>1$. Then there exist homogeneous polynomials $g_1,g_r$
of degrees $d_1,d_r$ which are polynomial combinations of the $f_i$
and have no common factor.
\endproclaim
\demo{Proof} Let $I=(f_i)_{1=2}^{i=r}$ be the homogeneous ideal generated by
$f_2,f_3, \dots, f_r$. Factor $f_1=\pi{p_i}$ where the $p_i$ are irreducible.
Let $[p_i] \subset {\p}{I_{d_r}}$ denote the linear subspace corresponding
to multiples of $p_i$. Then $[p_i]$ is a proper subspace, as otherwise
the $f_i$ would cut out a subscheme containing $Z(p_i)$, which has codimension
one. It follows that $\cup[p_i]$ is also a proper subspace, and hence
there exists $x \in {\p}{I_{d_r}}$ which avoids it. This $x$ gives rise
to a homogeneous polynomial $g_r$ of degree $d_r$ which has no common
factor with $g_1=f_1$.
\enddemo
\proclaim{Theorem 3.31} Let $\L$ be an even linkage class corresponding to the
stable equivalence class $[\N_0]$ via $\N$-type resolution, with
$\N_0 \neq 0$. Let $X_0$
be a minimal subscheme for $\L$ having $\N$-type resolution
$$0 \ra \oplus\O(-n)^{q(n)} \ra \N_0 \ra \I_{X_0}(h) \ra 0.$$
Then there exist hypersurfaces of degrees $a_0+h$ and $a_1+h$ which
contain $X_0$ and link $X_0$ to a subscheme $Y_0$ which is minimal
for its even linkage class.
\endproclaim
\demo{Proof} By definition of $a_1$, the morphism $\sigma_{a_1}$ is
surjective on an open set which contains the generic points of each
integral divisor in $\Pn$. It follows that this is also true for
the induced map $\oplus_{n \leq a_1}\O(-n)^{l(n)-q(n)} \ra \I_{X_0}(h)$.
Since the minimal value of $n$ making $l(n)-q(n) \neq 0$ is $n=a_0$,
it follows that this yields hypersurfaces of degrees varying from
$a_0+h$ to $a_1+h$ which cut out a subscheme of codimension $>1$.
Applying the previous lemma, we find that there are hypersurfaces of
degrees $a_0+h$ and $a_1+h$ which meet properly, and hence link $X_0$
to another subscheme $Y_0$. \par
Applying proposition 1.8 we obtain a resolution for $Y_0$ of
the form
$$0 \ra \oplus\O(n)^{l(n)-q(n)} \ra \E_0^\vee \oplus \O(-a_0) \oplus \O(-a_1)
\ra \I_{Y_0}(a_0+a_1+h) \ra 0$$
Proposition 3.29 shows that if $\M=\E_0^\vee \oplus \O(-a_0) \oplus \O(-a_1)$,
then $q_\M(n)=l(n)-q(n)$ for each $n \in \Z$. It follows from this
resolution and theorem 3.26 that $Y_0$ is minimal for its even
linkage class.
\enddemo
\heading References \endheading

\ref
\no 1 \by E. Ballico, G. Bolondi and J. Migliore \pages 117-128
\paper The Lazarsfeld-Rao problem for liaison classes of
two-codimensional subschemes of $\Pn$
\yr1991 \vol 113 \jour Amer J. Math.
\endref

\ref
\no 2 \by G. Bolondi and J. Migliore \pages 1-37
\paper The Structure of an Even Liaison Class
\yr 1989 \vol 316 \jour Trans. AMS
\endref

\ref
\no 3 \by A. Grothendieck
\paper El\'ements de G\'eome\`etrie Alg\'ebrique IV
\jour I.H.E.S. Publ. Math. \vol 28 \yr 1966
\endref

\ref
\no 4 \by R. Hartshorne \jour SLN \vol 156 \yr 1970
\paper Ample Subvarieties of Algebraic Varieties
\endref
\ref
\no 5 \by R. Hartshorne \book Algebraic Geometry \publ Springer-Verlag
\publaddr Berlin, Heidelberg and New York \yr 1977
\endref
\ref
\no 6 \by R. Hartshorne \paper Stable Reflexive Sheaves \yr 1980
\pages 121-176 \vol 254 \jour Math. Ann.
\endref
\ref
\no 7 \by R. Hartshorne \paper Generalized Divisors on Gorenstein Schemes
\jour K-theory \vol 8 \yr 1994 \pages 287-339
\endref
\ref
\no 8 \by S. Kleiman \paper Geometry on Grassmanians and applications to
splitting bundles and smoothing cycles \jour I.H.E.S. Publ. Math.
\vol 36 \yr 1969 \pages 281-298
\endref
\ref
\no 9 \by R. Lazarsfeld and A. P. Rao \yr 1983 \vol 997
\paper Linkage of General Curves of Large Degree \jour SLN \pages 267-289
\endref
\ref
\no 10 \by M. Martin-Deschamps and D. Perrin \publ Ast\'erisque \yr 1990
\book Sur la Classification des Courbes Gauches
\endref
\ref \no 11 \by M. Martin-Deschamps and D. Perrin \publ LMENS \vol 22
\yr 1992 \book Construction de Courves Lisses: un Th\'eor\`eme la Bertini
\endref
\ref \no 12 \by J. Migliore \yr 1994
\paper An Introduction to Deficiency Modules and Liaison Theory for
Subschemes of Projective Space
\endref
\ref \no 13 \by C. Peskine and L. Szpiro \jour Inventiones Math. \vol 26
\yr 1972 \pages 271-302 \paper Liaison des vari\'et\'es alg\'ebriques, I
\endref
\ref \no 14 \by A. P. Rao \paper Liaison Anong Curves in $\Pthree$
\jour Inventiones Math. \vol 50 \yr 1979 \pages 205-217
\endref
\ref \no 15 \by A. P. Rao \paper Liaison Equivalence Classes
\jour Math. Ann. \vol 258 \yr 1981 \pages 169-173
\endref
\enddocument